\documentclass[12pt]{amsart}

\usepackage{graphicx,url,lineno}
\usepackage{amssymb}
\usepackage{epstopdf}
\usepackage{color}
\usepackage{mathrsfs}

\newtheorem{theorem}{Theorem}
\newtheorem{lemma}{Lemma}
\newtheorem{definition}{Definition}
\newtheorem{remark}{Remark}

\newtheorem{corollary}{Corollary}
\newtheorem{proposition}{Proposition}

\title[Periodic functions and fractional order operators]{On quasi-periodicity properties of fractional integrals and fractional derivatives of periodic functions}

\author[Area]{I. Area}
\author[Losada]{J. Losada}
\author[Nieto]{Juan J. Nieto}

\address[Area]{Departamento de Matem\'{a}tica Aplicada II, E.E. Telecomunicaci\'{o}n, Universidade de Vigo, 36310-Vigo, Spain.}\email[Area]{area@uvigo.es}

\address[Losada]{Facultade de Matem\'{a}ticas, Universidade de Santiago de
Compostela, 15782-Santiago de Compostela, Spain. Corresponding author}
\email[Losada]{jorge.losada@rai.usc.es}

\address[Nieto]{Facultade de Matem\'{a}ticas, Universidade de Santiago de
Compostela, 15782-Santiago de Compostela, Spain, and Faculty of Science, King Abdulaziz University, P.O. Box 80203, 21589, Jeddah, Saudi Arabia.}
\email[Nieto]{juanjose.nieto.roig@usc.es}

\date{\today}                                        

\begin{document}

\begin{abstract}
This paper is devoted to the study of quasi-periodic properties of fractional order integrals and derivatives of periodic functions. Considering Riemann-Liouville and Caputo definitions, we discuss when the fractional derivative and when the fractional integral of a certain class of periodic functions satisfies particular properties. We study concepts close to the well known idea of periodic function, such as S-asymptotically periodic, asymptotically periodic or almost periodic function. Boundedness of fractional derivative and fractional integral of a periodic function is also studied.
\end{abstract}

\maketitle

\section{Introduction}

Periodic functions play a central role in mathematics since the seminal works of J.\,B. Fourier. Indeed, the study of existence of periodic solutions is one of the most interesting and important topics in qualitative theory of differential equations. This is due to its implications in pure or abstract areas of mathematics, but also due to its applications, ranging from physics to natural and social sciences and of course, in control theory.

However, the definition of periodic function is extremely demanding and then, the conditions to guarantee the existence of periodic solutions are very harsh. For this reason, in the past decades, many authors (see \cite{henriquez2,henriquez1, nicola-pierri} and references therein) have proposed and studied extensions of the concept of periodicity which have shown interesting and useful.

The idea of integrals and derivatives of noninteger order goes back to the birth of the theory of differential calculus, more exactly to the Leibniz's note in his letter to L'H\^{o}pital dated 30 September 1695 \cite{oldham}. Probably, the first application of fractional calculus was made by N.\,H. Abel in 1823 when he was studying the integral equation that arises in the formulation of the tautochrone problem. But we also like to mention O. Heaviside, who in 1893 used fractional differential operators to study the {\it{Age of the Earth}} (see \cite{terra} or \cite[Chapter 7]{duarte}).

For a long time, the theory of Fractional Calculus developed only as a pure theoretical field of mathematics. However, in the last decades, it was found that fractional derivatives and integrals provide, in some situations, a better tool to understand some physical phenomena, especially when dealing with processes with memory \cite{memory}. Applications includes modeling viscoelastic and viscoplastic materials \cite{death}, chemical processes \cite{brusselator}, and a wide range of engineering problems. Approximations of fractional derivatives have been studied in \cite{baeumer}, fractional Sobolev spaces in \cite{Schikorra}, and the fractional Laplacian in \cite{cabresire,seoi,servadei}.

It is an obvious fact that the classical derivative, if it exists, of a periodic function is also a periodic function of the same period. Also the primitive of a periodic function may be periodic (e.g., $-\cos t$ as primitive of $\sin t$). Nevertheless, when we consider derivatives or integrals of non integer order, this fact is not true \cite{ALOSN}. Periodicity is also important in the context of fractional calculus \cite{ALOSN,MR2863974,khanlamb,love,maxprinc}.

Our aim in this paper is to study quasi-periodic properties of fractional order integrals and derivatives of periodic functions. Using Riemann-Liouville and Caputo definitions, we take into consideration concepts as asymptotically periodic function, S-asymptotically periodic function or almost periodic function. 

This work has six sections. In the next section we introduce the basic definitions. In Section \ref{S:3}, we present some useful results that will be used in Section \ref{S:4}, where we prove our main conclusions. In Section \ref{S:5} we comment for a classical example of periodic function the results previously proved. Concluding remarks in Section \ref{S:7} close the paper.

\section{Preliminaries}

The fractional integral of order $\alpha>0$ of a given function $f\colon\mathbb{R}\longrightarrow\mathbb{R}$ is defined as \cite{MR2218073,MR1658022}
\begin{equation}\label{eq:1}
I^\alpha f\,(t)=\dfrac{1}{\Gamma(\alpha)}\int_0^t(t-s)^{\alpha-1}f(s)\,ds,\qquad t>0.
\end{equation}
Note that for $\alpha<1$ the integral may be singular, but it is well defined if, for example,  $f \in {{\rm L}}^{1}_{\text{loc}}({\mathbb{R}})$. Using (\ref{eq:1}), the fractional Riemann-Liouville derivative of order $\alpha$, $\alpha\in(0,1)$, of $f$ is defined as (see \cite{MR2218073,MR1658022})
\[{}^{\textsc{rl}}D^{\alpha}f\,(t)=D^{1}\left(I^{1-\alpha} f\right)\,(t)=\frac{1}{\Gamma(1-\alpha)}\frac{d}{dt} \int_{0}^{t} (t-s)^{-\alpha} f(s)ds,\]
provided the right hand side is defined for almost every $t\in\mathbb{R}^+$. This is well defined if, for example, $f$ is absolutely continuous on every compact interval of $\mathbb{R}$.

There are many more definitions of fractional integral and fractional derivative. We are not giving a complete list, but recall the Caputo fractional derivative \cite{MR2218073,MR1658022}
\begin{equation*}
{}^{\textsc c}D^{\alpha}f\,(t)={}^{\textsc{rl}}D^{\alpha}g\,(t),
\end{equation*}
with $g(t)=f(t)-f(0)$.
In addition, if $f$ is an absolutely continuous function on every compact interval (of $\mathbb{R}$), we can write, for $\alpha\in(0,1)$,
\begin{equation}\label{eq:3}
{}^{\textsc c}D^{\alpha}f\,(t)= I^{1-\alpha}D^{1}f\,(t)=\dfrac{1}{\Gamma(1-\alpha)}\int_0^t(t-s)^{-\alpha}f'(s)\,ds.
\end{equation}

Note also that if $\alpha\in (0,1)$ and $f$ is a function for which the Caputo fractional derivative, ${}^{\textsc c}D^\alpha f$, exists together with the Riemann-Liouville fractional derivative, ${}^{\textsc{rl}}D^\alpha f$, both of them of order $\alpha$, then we have that \cite[(2.4.8), p. 91]{MR2218073}
\[{}^{\textsc c}D^\alpha f\,(t)={}^{\textsc{rl}}D^\alpha f\,(t)-\dfrac{f(0)}{\Gamma(1-\alpha)}\,t^\alpha.\]
Moreover, as in the integer case, we have
\[^{\textsc {rl}}D^{\alpha} (I^{\alpha}f)\,(t)=f(t), \qquad \,{}^{\textsc{c}}D^{\alpha}( I^{\alpha}f )\,(t)=f(t);\]
but $I^{\alpha}(^{\textsc{rl}}D^{\alpha}f)$ or $I^{\alpha}({}^{\textsc{c}}D^{\alpha}f)$ are not, in general, equal to $f$.

It will be also useful the definition of Weyl fractional integral. Following \cite[Chapter 4, Section 19]{samko1993fractional}, the Weyl fractional integral of a function $f$, of order $\alpha\in (0,1)$ is given by
\[^{\textsc w}I^\alpha f\,(t)=\dfrac{1}{2\pi}\int_0^{2\pi}f(t-s)g(s)\,ds,\]
where, for $0\le s\le 2\pi$,
\[g(s)=\dfrac{2\pi}{\Gamma(\alpha)}s^{\alpha-1}+\dfrac{1}{\Gamma(\alpha)}\lim_{n\to \infty}\left[2\pi\sum_{m=1}^n(s+2\pi m)^{\alpha-1}-\dfrac{(2\pi n)^\alpha}{\alpha}\right].\]
We recall that (see \cite[Lemma 19.3]{samko1993fractional}) if $f\colon\mathbb{R}\longrightarrow\mathbb{R}$ is a $T$-periodic function and $f \in{\rm L}^1(0,T)$ is such that
\begin{equation}\label{eq:meanzero}
\int_0^T f(t)\,dt=0,
\end{equation}
then, for $\alpha\in (0,1)$,
\begin{equation}\label{eq:4}
^{\textsc w}I^\alpha f\,(t)=\dfrac{1}{\Gamma(\alpha)}\int_{-\infty}^t (t-s)^{\alpha-1}f(s)\,ds,
\end{equation}
provided that the right-hand side integral is understood as conventionally convergent:
\[\int_{-\infty}^t (t-s)^{\alpha-1}f(s)\,ds=\lim_{\substack{n\to\infty\\ n\in\mathbb{N}}}\int_{t-nT}^t (t-s)^{\alpha-1}f(s)\,ds.\]\medskip

Following now \cite{nicola-pierri}, we continue recalling some useful concepts, definitions and results about quasi-periodic functions. Denote by $\mathbb{R}^+$ the infinite interval $[0,+\infty)$, $\mathcal{C}\left(\mathbb{R},\mathbb{R}\right)$ is the space of all continuous functions defined on $\mathbb{R}$, $\mathcal{C}_b\left(\mathbb{R}^+,\mathbb{R}\right)$ is the space of all continuous and bounded functions from $\mathbb{R}^+$ into $\mathbb{R}$ endowed with the norm of the uniform convergence denoted by $\|\cdot\|_\infty$. Its subspaces, $\mathcal{C}_0\left(\mathbb{R}^+,\mathbb{R}\right)$ and $\mathcal{C}_T\left(\mathbb{R}^+,\mathbb{R}\right)$, are defined as
\begin{align*}
\mathcal{C}_0\left(\mathbb{R}^+,\mathbb{R}\right)&=\left\{f\in\mathcal{C}_b\left(\mathbb{R}^+,\mathbb{R}\right)\,\colon\, \lim_{t\to\infty}|f(t)|=0\right\},\\
\mathcal{C}_T\left(\mathbb{R}^+,\mathbb{R}\right)&=\left\{f\in\mathcal{C}_b\left(\mathbb{R}^+,\mathbb{R}\right)\,\colon\, f\text{ is $T$-periodic}\right\}.
\end{align*}
We denote by $\mathcal{C}_T(\mathbb{R},\mathbb{R})$ the linear space of all continuous and $T$-periodic functions.

\begin{definition}\cite{besicovitch}
A set $E\subset\mathbb{R}$ is said relatively dense if it exists a number $l>0$ such that any interval of length $l$ contains at least on number of $E$.
\end{definition}

\begin{definition}\cite{besicovitch}
A function $f\in\mathcal{C}\left(\mathbb{R},\mathbb{R}\right)$ is almost periodic if for every $\varepsilon>0$ there exists a relatively dense subset $\mathscr{H}(\varepsilon,f)$ of\, $\mathbb{R}$ such that $\left|f(t+\xi)-f(t)\right|<\varepsilon$, for every $t\in\mathbb{R}$ and all $\xi\in\mathscr{H}(\varepsilon,f)$. \\
We denote by $AP(\mathbb{R})$ the space consisting of all almost periodic functions.
\end{definition}

\begin{definition}
A function $f\in\mathcal{C}\left(\mathbb{R},\mathbb{R}\right)$ is called asymptotically almost periodic if there exists an almost periodic function $f_1$ and a function $f_2\in\mathcal{C}_0\left(\mathbb{R}^+,\mathbb{R}\right)$ such that $f=f_1+f_2$. If $f_1$ is periodic (resp. $T$-periodic) $f$ is said to be asymptotically periodic (resp. asymptotically $T$-periodic).\\
We denote by $AAP(\mathbb{R})$ the space of all asymptotically periodic functions  and by $AP_T(\mathbb{R})$ the space of asymptotically $T$-periodic functions.
\end{definition}

\begin{definition}
A function $f\in\mathcal{C}_b\left(\mathbb{R}^+,\mathbb{R}\right)$ is called S-asymptotically $T$-periodic if there exists $T>0$ such that
\[\lim_{t\to+\infty}\left[f(t+T)-f(t)\right]=0.\]
In this case we say that $T$ is an asymptotically period of $f$.\\
We denote by $SAP_{T}(\mathbb{R})$ the space of all S-asymptotically $T$-periodic functions.
\end{definition}
In \cite{henriquez1,nicola-pierri} authors show examples of functions that are S-asymptotically periodic but not asymptotically periodic.

\section{Some useful results}\label{S:3}

To obtain our results about quasi-periodicity properties of fractional derivatives and integrals, we need some previous results. In this section we recall and present some facts that will be used in next section.

\begin{theorem}\cite[Theorem 1]{ALOSN}\label{T:noperiodic}
Let $f\in\mathcal{C}_T(\mathbb{R},\mathbb{R})$ a nonzero $T$-periodic function with $f\in{\rm L}^1_{\text{\rm loc}}(\mathbb{R})$. Then $I^\alpha f$ cannot be a $T$-periodic function.
\end{theorem}

The above result has been proved for the first time by Tavazoei \cite{Tavazoei}, Kaslik and Sivasundaram have also proved the same result by using the Mellin transform in \cite{MR2863974}. The same result appears in \cite{MR2971825} using the Laplace transform. Recently \cite{ALOSN} we have proved that $I^{\alpha} f$ cannot be periodic for any period $\tilde{T}$ unless, of course, that $f\equiv 0$.\medskip

\begin{lemma}\label{L:2.5}
If $f\in \mathcal{C}_T(\mathbb{R},\mathbb{R})$ and $\alpha\in (0,1)$, we have for $t>0$ that
\[\lim_{\substack{n\to\infty\\n\in\mathbb{N}}}\int_{t-nT}^{t}(t-s)^{\alpha-1}f(s)\,ds=\lim_{\substack{n\to\infty\\n\in\mathbb{N}}}\int_{-nT}^{t}(t-s)^{\alpha-1}f(s)\,ds.\]
\begin{proof}
It is clear that
\begin{equation*}
\left|\int_{-nT}^{t-nT}(t-s)^{\alpha-1}f(s)\,ds\right|\le \int_{-nT}^{t-nT}(t-s)^{\alpha-1}|f(s)|\,ds\le t\,(nT)^{\alpha-1}\,\|f\|_{\infty};
\end{equation*}
so, since $\alpha\in (0,1)$, we deduce that for all $t>0$
\[\lim_{\substack{n\to\infty\\n\in\mathbb{N}}} \int_{-nT}^{t-nT}(t-s)^{\alpha-1}f(s)\,ds =0.\]
The proof follows taking into account that
\[\int_{t-nT}^{t}(t-s)^{\alpha-1}f(s)\,ds=\int_{-nT}^{t}(t-s)^{\alpha-1}f(s)\,ds-\int_{-nT}^{t-nT}(t-s)^{\alpha-1}f(s)\,ds.\qedhere\]
\end{proof}
\end{lemma}

Next three results will be useful to deduce that the fractional integral of a periodic function is not an almost periodic function.

\begin{lemma}\cite[Lemma 3.1]{henriquez1}\label{L:qp1}
Let $f\colon\mathbb{R}\longrightarrow \mathbb{R}$ be an S-asymptotically $T$-periodic function, and $\left(t_n\right)_{n\in\mathbb{N}}$ be a sequence with $t_n\to+\infty$ as $n\to+\infty$ and assume that $f_{t_n}(t)=f(t+t_n)$ satisfies that $f_{t_n}\to F$ uniformly on compact subsets of\, $\mathbb{R}^+$. Then $F\in\mathcal{C}_T(\mathbb{R},\mathbb{R})$.
\begin{proof}
It is clear that $F$ is a continuous function. For $T\ge 0$ and $\varepsilon>0$, we select $n_0\in\mathbb{N}$ such that
\begin{align*}
&\left|F(s)-f(s+t_n)\right|\le\varepsilon,\qquad s\in[t,t+T],\\
&\left|f(t+t_n)-F(t)\right|\le\varepsilon,\qquad t \ge 0,
\end{align*}
for every $n\ge n_0$. Hence, for $n\ge n_0$, we have that
\begin{align*}
\left|F(t+T)-F(t)\right|\le& \left|F(t+T)-f(t+T+t_n)\right|+\left|f(t+T+t_n)-f(t+t_n)\right|\\
&+\left|f(t+t_n)-F(t)\right|\le 3\,\varepsilon,
\end{align*}
which implies that $F(t+T)=F(t)$. The proof is complete.
\end{proof}
\end{lemma}

\begin{proposition}\cite[Proposition 3.4]{henriquez1}\label{L:qp2}
Let $f\colon\mathbb{R}^+\longrightarrow\mathbb{R}$ be an S-asymptotically $T$-periodic and asymptotically almost periodic function. Then $f$ is an asymptotically $T$-periodic function.
\begin{proof}
We can decompose $f=f_1+f_2$, where $f_1$ is an almost periodic function and $f_2\in\mathcal{C}_0(\mathbb{R}^+,\mathbb{R})$. It follows from the theory of almost periodic functions that there exists a sequence of real numbers $\left(t_n\right)_{n\in\mathbb{N}}$ such that $t_n\to+\infty$ and ${f_1}_{t_n}(t)=f_1(t+t_n)$ satisfies that ${f_1}_{t_n}\to f_1$ as $n\to+\infty$ uniformly on $\mathbb{R}^+$.

Therefore ${f_1}_{t_n}\to f_1$ as $n\to +\infty$ uniformly on $\mathbb{R}^+$ and, by Lemma \ref{L:qp1}, it follows that $f_1\in\mathcal{C}_T(\mathbb{R},\mathbb{R})$ which in turns implies that the function $f$ is asymptotically $T$-periodic.
\end{proof}
\end{proposition}

\begin{proposition}\label{L:qp3}
Let $f\colon\mathbb{R}\longrightarrow\mathbb{R}$ be an S-asymptotically $T$-periodic and almost periodic function. Then $f$ is a $T$-periodic function.
\begin{proof}
As previously, it follows from the theory of almost periodic functions that there exists a sequence of real numbers $\left(t_n\right)_{n\in\mathbb{N}}$ such that $t_n\to+\infty$ and $f_{t_n}\to f$ as $n\to+\infty$ uniformly on $\mathbb{R^+}$. Applying Lemma \ref{L:qp1}, it follows that $f\in\mathcal{C}_T(\mathbb{R},\mathbb{R})$.
\end{proof}
\end{proposition}

\section{Main results}\label{S:4}

The main purpose of this paper is to study if the fractional integral (or derivative) of a periodic function possesses a {\it{quasi-periodic}} behavior. We will answer the following relevant questions: if $f\in\mathcal{C}_T(\mathbb{R},\mathbb{R})$,
\begin{itemize}
\item[1.-] Is then $I^{\alpha} f$ bounded?
\item[2.-] Is then $I^{\alpha} f$ S-asymptotically $T$-periodic?
\item[3.-] Is then $I^{\alpha} f$ asymptotically $T$-periodic?
\item[4.-] Is then $I^{\alpha} f$ almost-periodic?
\end{itemize}
We also present some results concerning fractional derivatives.\medbreak

We begin this section with a result about S-asymptotically periodic functions and fractional integrals. Then, we focus our attention in the concept of asymptotically periodicity

\begin{lemma}\label{L:1}
If $f\in\mathcal{C}_T\left(\mathbb{R}^+,\mathbb{R}\right)$ and $\alpha\in(0,1)$ then $I^\alpha f$ satisfies the following property of asymptotic $T$-periodicity:
\begin{equation}
\nonumber
\lim_{t \to +\infty}\left[I^\alpha f\,(t+T)-I^\alpha f\,(t)\right]=0.
\end{equation}
\begin{proof}
Indeed, for $t \geq 0$, since $f$ is a $T$-periodic function, we have that
\begin{align*}
I^\alpha f&\,(t+T)-I^\alpha f\,(t)=\\
&=\frac{1}{\Gamma(\alpha)}\left(\int_{0}^{t+T}(t+T-s)^{\alpha-1}f(s)\,ds-\int_{0}^{t}(t-s)^{\alpha-1}f(s)\,ds\right)\\
&=\frac{1}{\Gamma(\alpha)}\left(\int_{0}^{t+T}(t+T-s)^{\alpha-1}f(s)\,ds-\int_{T}^{t+T}(t-r+T)^{\alpha-1}f(r-T)\,dr\right)\\
&=\frac{1}{\Gamma(\alpha)}\left(\int_{0}^{t+T} (t+T-s)^{\alpha-1}f(s)\,ds-\int_{T}^{t+T} (t-s+T)^{\alpha-1}f(s)\,ds\right).
\end{align*}
Thus,
\begin{equation*}
I^\alpha f\,(t+T)-I^\alpha f\,(t)=\frac{1}{\Gamma(\alpha)}\int_{0}^{T}(t+T-s)^{\alpha-1}f(s)\,ds,
\end{equation*}
which implies that
\begin{equation*}
\left|I^\alpha f\,(t+T)-I^\alpha f\,(t)\right|\leq\frac{1}{\Gamma(\alpha)}\,T\,\|f\|_{\infty}t^{\alpha-1}\le T\,\|f\|_{\infty}t^{\alpha-1},
\end{equation*}
since $0<1/\Gamma(\alpha)<1$ for $0<\alpha<1$. The conclusion of this Lemma follows now easily.
\end{proof}
\end{lemma}

First of all, we would like to notice here that this result is contained in \cite{MR2971825}.
\begin{remark}\label{R:1}
Note also that we have not proved that $I^\alpha f$ is an S-asymptotically $T$-periodic function, because $I^\alpha f$ may be an unbounded function.
\end{remark}

\begin{theorem}\label{T:4.2}
Let $f\in\mathcal{C}_T\left(\mathbb{R}^+,\mathbb{R}\right)$ such that $I^\alpha f$, with $\alpha\in(0,1)$, is a bounded function. Then $I^\alpha f$ is an S-asymptotically $T$-periodic function.
\begin{proof}
Obvious from Lemma \ref{L:1} and Remark \ref{R:1}.
\end{proof}
\end{theorem}

\begin{lemma}\label{C:1}
Let $f\in\mathcal{C}_T\left(\mathbb{R}^+,\mathbb{R}\right)$ an absolutely continuous function. Then for $\alpha\in (0,1)$, ${}^{\textsc c}D^\alpha f$ satisfies the following asymptotic $T$-periodic property,
\[\lim_{t\to\infty}\left({}^{\textsc c}D^\alpha f\,(t+T)-{}^{\textsc c}D^\alpha f\,(t)\right)=0.\]
\begin{proof}
It is clear from the definition of Caputo fractional derivative, (\ref{eq:3}) and Lemma \ref{L:1} since the ordinary derivative of a $T$-periodic function also is a $T$-periodic function.
\end{proof}
\end{lemma}

\begin{remark}\label{R:2}
We have not proved that ${}^{\textsc{c}}D^\alpha f$ is an S-asymptotically $T$-periodic function, because we do not know if ${^\textsc{c}}D^\alpha f$ is bounded.
\end{remark}

\begin{theorem}
Let $f\in\mathcal{C}_T\left(\mathbb{R}^+,\mathbb{R}\right)$ an absolutely continuous function and $\alpha\in (0,1)$. Then, if ${}^{\textsc c}D^\alpha f$ is a bounded function, it is an S-asymptotically $T$-periodic function.
\begin{proof}
Obvious from Lemma \ref{C:1} and Remark \ref{R:2}.
\end{proof}
\end{theorem}

\begin{remark}
If we suppose in addition that $f$ is such that $f^{(n-1)}$ is an absolutely continuous function ($f\in AC^{n}$) then, if $^{\textsc c}D^\alpha f$ with $\alpha\in (n-1,n)$ is bounded, we have that $^{\textsc c}D^\alpha f$ is an S-asymptotically $T$-periodic function.
\end{remark}

Now, we will prove a similar result for Riemann-Liouville fractional derivative.

\begin{theorem}
Let $f$ a function such that $f\in \mathcal{C}_T(\mathbb{R},\mathbb{R})$ and $f\in AC^{n}(\mathbb{R})$. Then, if\, $^{\textsc{rl}}D^\alpha f$ with $\alpha\in(n-1,n)$ is a bounded function, we have that $^{\textsc{rl}}D^\alpha f$ is an S-asymptotically $T$-periodic function.
\begin{proof}
Suppose that there exists a function $f$ in the conditions above, but such that $^{\textsc{rl}}D^\alpha f$ is not an  S-asymptotically $T$-periodic function. Consider in that case the function $\tilde{f}=f-f(0)$; it is obvious that $\tilde{f}\in \mathcal{C}_T(\mathbb{R},\mathbb{R})$ and $\tilde{f}\in AC^{n}(\mathbb{R})$. Taking into account that
\[^{\textsc{c}}D^\alpha \tilde{f}\,(t)={}^{\textsc{rl}}D^\alpha f\,(t),\]
we obtain that $^{\textsc{c}}D^\alpha \tilde{f}$ is not an S-asymptotically $T$-periodic function, which is not true.
\end{proof}
\end{theorem}\medskip

In the next theorem, we will assume that $I^\alpha f$ is a bounded function.

\begin{theorem}\label{T:1}
Let $f\in\mathcal{C}_T(\mathbb{R},\mathbb{R})$ such that $I^\alpha f$, with $\alpha\in(0,1)$, is a bounded function. Then, $I^\alpha f$ is an asymptotically $T$-periodic function.

\begin{proof}
(Following \cite[Lemma 2.1]{haiyin}) Let for brevity
\[\varphi(t)=I^\alpha f\,(t)=\dfrac{1}{\Gamma(\alpha)}\int_0^t(t-s)^{\alpha-1} f(s)\,ds,\qquad t>0.\]

For each $n\in\mathbb{N}$, we consider the following functions: $\varphi_n(t)=\varphi(t+nT)$ and $\Phi_n(t)=\sup_{k\ge n}\varphi_k(t)$ for $t\ge 0$. Since $I^\alpha f$ is a bounded function, $\varphi_n$ an $\Phi_n$ are bounded and continuous functions  defined on $\mathbb{R}^+$. Moreover, we have that $\Phi_{n+1}(t)\le \Phi_n(t)$ for all $n\in\mathbb{N}$ and $t\in\mathbb{R}^+$. In addition, from the hypothesis of $T$-periodicity over $f$, we deduce that
\begin{align*}
\lim_{n\to\infty}\Phi_n(t)&=\lim_{n\to \infty}\left[\sup_{k\ge n} \varphi(t+kT)\right]=\limsup_{n\to \infty}\varphi(t+nT)\\
&=\limsup_{n\to \infty}\dfrac{1}{\Gamma(\alpha)}\int_0^{t+nT}(t+nT-s)^{\alpha-1}f(s)\,ds\\
&=\limsup_{n\to \infty}\dfrac{1}{\Gamma(\alpha)}\int_{-nT}^{t}(t-r)^{\alpha-1}f(r+nT)\,dr\\
&=\limsup_{n\to \infty}\dfrac{1}{\Gamma(\alpha)}\int_{t-nT}^{t}(t-s)^{\alpha-1}f(s)\,ds,
\end{align*}
where we have used Lemma \ref{L:2.5}.

Taking now into account that for each $n\in\mathbb{N}$,
\begin{align*}
\Big|\int_{t-nT}^t(t-s)^{\alpha-1}f(s)\,ds&-\int_{t-(n+1)T}^t(t-s)^{\alpha-1}f(s)\,ds\Big|\\
&=\Big|\int_{t-(n+1)T}^{t-nT}(t-s)^{\alpha-1}f(s)\,ds\Big|\le T\,(nT)^{\alpha-1}\,\|f\|_{\infty},
\end{align*}
we have that, for each $t\in\mathbb{R}^+$, the sequence $(a_n)_{n\in\mathbb{N}}$ given for each $n\in\mathbb{N}$ by
\[a_n=\int_{t-nT}^t (t-s)^{\alpha-1}f(s)\,ds,\]
is a Cauchy sequence in $\left(\mathbb{R},|\cdot|\right)$. Therefore,
\[\lim_{n\to\infty}\Phi_n(t)=\Phi(t)\in \mathbb{R},\quad\text{for all }t\ge 0\]
and by (\ref{eq:4}), we obtain that
\begin{equation}\label{eq:liouville}
\Phi(t)\equiv {}^{\textsc w}I^\alpha f\,(t)=\dfrac{1}{\Gamma(\alpha)}\int_{-\infty}^t(t-s)^{\alpha-1}f(s)\,ds,\qquad t\ge 0,
\end{equation}
is the Weyl's integral (sometimes also called Liouville fractional integral) of order $\alpha$ of $f$. Note that we are only considering $f\in\mathcal{C}_T(\mathbb{R},\mathbb{R})$, but not necessarily with mean zero.

In other words, we have that $\left(\Phi_n\right)_{n\in\mathbb{N}}$ is a sequence of functions that converges pointwise to the continuous function given by (\ref{eq:liouville}).

In such case, all hypothesis of Dini's Theorem (see \cite{bartle,rudin1964principles}) are satisfied and we deduce that $\left(\Phi_n\right)_{n\in\mathbb{N}}$ uniformly converges to $\Phi$ on any given closed subinterval of\, $\mathbb{R}^+$.\medbreak

We show now that $\Phi$ is a $T$-periodic function. By Lemma \ref{L:1} and Heine Theorem, we have that
\[\lim_{n\to\infty}\left(\varphi(t+nT+T)-\varphi(t+nT)\right)=0\]
and thus
\begin{align*}
\Phi(t+T)-\Phi(t)&=\lim_{n\to\infty}\left(\Phi_n(t+T)-\Phi_n(t)\right)\\
&=\lim_{n\to\infty}\left(\sup_{k\ge n}\varphi(t+T+kT)-\sup_{k\ge n}\varphi(t+kT)\right)\\
&=\limsup_{n\to\infty}\left(\varphi(t+T+nT)-\varphi(t+nT)\right)\\
&=\lim_{n\to\infty}\left(\varphi(t+T+nT)-\varphi(t+nT)\right)=0.
\end{align*}

Therefore, we have that $\Phi$ is a $T$-periodic function. This fact seems to be already known when the mean of $f$ is zero (see, for example, \cite[p. 348]{samko1993fractional}), that is, when (\ref{eq:meanzero}) holds.\medbreak

Next step, we show that $(\Phi_n)_{n\in\mathbb{N}}$ uniformly converges to $\Phi$ on $\mathbb{R}^+$.

Let for each $m,\,n\in\mathbb{N}$, $I_{n}=\left[n\,T,(n+1)T\right]$ and
\[A_{m,\,n}=\sup_{t\in I_m}\left|\Phi_n(t)-\Phi(t)\right|.\]
Using now that $\Phi$ is a $T$-periodic function, we have that
\begin{align*}
A_{m,\,n}&=\sup_{t\in I_m}\left|\sup_{k\ge n}\varphi_k(t)-\Phi(t)\right|=\sup_{t\in I_m}\left|\sup_{k\ge n}\varphi(t+kT)-\Phi(t)\right|\\
&=\sup_{t\in I_0}\left|\sup_{k\ge n+m}\varphi(t+kT)-\Phi(t)\right|=A_{0,\,m+n}.
\end{align*}
Since $I^\alpha f$ is bounded, it is clear that $\lim_{s\to\infty}A_{0,\,s}=0$. Therefore, for any given $\varepsilon>0$, there exists $N\in\mathbb{N}$ such that for any $k\ge N$
\begin{equation}\label{eq:6}
A_{0,\,k}=\sup_{t\in I_0}\left|\Phi_k(t)-\Phi(t)\right|<\varepsilon.
\end{equation}
Given $t\in\mathbb{R}^+$ arbitrary, there exists $m\in\mathbb{N}$ such that $t\in I_m$ and from (\ref{eq:6}), when $n\ge N$,
\[A_{m,\,n}=A_{0,\,m+n}=\sup_{t\in I_{m}}\left|\Phi_n(t)-\Phi(t)\right|<\varepsilon.\]
That is to say, we have that for $n\ge N$
\[\sup_{t\in\mathbb{R}^+}\left|\Phi_n(t)-\Phi(t)\right|<\varepsilon,\]
or equivalently, $\left(\Phi_n\right)_{n\in\mathbb{N}}$ uniformly converges to $\Phi$ on $\mathbb{R}^+$.\medskip

If instead of $\Phi_n$ we consider for each $n\in\mathbb{N}$
\[\Psi_n(t)=\inf_{k\ge n}\varphi_k(t)=\inf_{k\ge n}\varphi(t+kT),\]
we can easily see, as before, that $\left(\Psi_n\right)_{n\in\mathbb{N}}$ also uniformly converges to $\Phi$ on $\mathbb{R}^+$.\medbreak

Next we will prove that
\[\lim_{t\to+\infty}\varphi(t)-\Phi(t)= 0.\]
Because of
\[\Psi_n(t)=\inf_{k\ge n}\varphi_k(t)\le\varphi_n(t)\le\sup_{k\ge n}\varphi_k(t)=\Phi_n(t)\]
and $\left(\Phi_n\right)_{n\in\mathbb{N}}$, $\left(\Psi_n\right)_{n\in\mathbb{N}}$ uniformly converges to $\Phi$ on $\mathbb{R}^+$ we conclude that $\left(\varphi_n\right)_{n\in\mathbb{N}}$ uniformly converges to $\Phi$ on $\mathbb{R}^+$.

Thus, for any given $\varepsilon>0$ there exists an $N\in\mathbb{N}$ such that for any $n\ge N$
\[\left|\varphi_n(t)-\Phi(t)\right|<\dfrac{\varepsilon}{2}\qquad\text{for all }t\in\mathbb{R}^+,\]
and we have
\begin{align}\label{eq:7}
\left|\varphi(t)-\Phi(t)\right|&\le\left|\varphi(t)-\varphi_n(t)\right|+\left|\varphi_n(t)-\Phi(t)\right|<\left|\varphi(t)-\varphi_n(t)\right|+\dfrac{\varepsilon}{2}.
\end{align}
{}From Lemma \ref{L:1}, for the above $\varepsilon>0$ there exists $S>0$ such that if $t>S$ then
\[\left|\varphi(t+T)-\varphi(t)\right|<\dfrac{\varepsilon}{2N}.\]
Using this, for $t>S$, it follows that
\begin{align}\label{eq:8}
\nonumber
\left|\varphi(t)-\varphi(t+nT)\right|&\le\left|\varphi(t)-\varphi(t+T)\right|+\left|\varphi(t+T)-\varphi(t+2T)\right|\\
\nonumber
&\hspace*{0.5cm}+\cdots+\left|\varphi(t+NT)-\varphi(t+(N-1)T)\right|\\
&<\dfrac{\varepsilon}{2 N}\,N=\dfrac{\varepsilon}{2}.
\end{align}

Therefore, taking into account (\ref{eq:7}) and (\ref{eq:8}), we conclude that, for $t>S$
\[\left|\varphi(t)-\Phi(t)\right|<\dfrac{\varepsilon}{2}+\dfrac{\varepsilon}{2}=\varepsilon\]
and then
\[\lim_{t\to\infty}\left(\Phi(t)-\varphi(t)\right)=0.\]\medskip

{}From the obvious identity $\varphi(t)=\Phi(t)+\left(\varphi(t)-\Phi(t)\right)$, we finally deduce that $\varphi$ is an asymptotically $T$-periodic function with $f_1=\Phi$ and $f_2=\varphi-\Phi$.
\end{proof}
\end{theorem}

\begin{theorem}
Let $\alpha>0$, $\alpha\in[n-1,n)$ with $n\in\mathbb{N}$. If $f\in\mathcal{C}_T(\mathbb{R},\mathbb{R})$, $f\in AC^n(\mathbb{R})$ and $I^{n-\alpha} f$ is a bounded function then ${}^{\textsc c}D^\alpha f$ is an asymptotically $T$-periodic function.
\begin{proof}
It is clear from the definition of Caputo fractional derivative and Theorem \ref{T:1}, since the ordinary derivative of a $T$-periodic function also is a $T$-periodic function.
\end{proof}
\end{theorem}

Next we will prove a similar result using Riemann-Liouville derivative.

\begin{theorem}
If $f\in\mathcal{C}_T(\mathbb{R},\mathbb{R})$ and $I^{1-\alpha}$, with $\alpha\in (0,1)$, is a bounded function, then ${}^{\textsc{rl}}D^\alpha f$ is an asymptotically $T$-periodic function.
\begin{proof}
Taking into account that ${}^{\textsc{rl}}D^\alpha f\,(t)=D^1\left(I^{1-\alpha}f\right)\,(t)$ and Theorem \ref{T:1}, we have
\[^{\textsc{rl}}D^\alpha f\,(t)=D^1\left(f_1+f_2\right)\,(t)={f_1}'(t)+{f_2}'(t),\]
where $f_1\in\mathcal{C}_T(\mathbb{R}^+,\mathbb{R})$ and $f_2\in\mathcal{C}_0(\mathbb{R}^+,\mathbb{R})$.

Since ${f_1}'\in\mathcal{C}_T(\mathbb{R}^+,\mathbb{R})$ and $^{\textsc{rl}}D^\alpha f$ is a bounded function, we have that
\[\max\left\{\limsup_{t\to +\infty}|{f_2}'(t)|,\, \liminf_{t\to +\infty}|{f_2}'(t)|\right\}<+\infty.\]
If
\[\limsup_{t\to +\infty}{f_2}'(t)=\liminf_{t\to +\infty}{f_2}'(t)=\lim_{t\to+\infty} {f_2}'(t)=d\in\mathbb{R},\]
then $^{\textsc{rl}}D^\alpha f\,(t)=({f_1}'(t)+d)+({f_2}'(t)-d)$, and the theorem is proved. In other case, that is, if
\[d_1=\limsup_{t\to +\infty}{f_2}'(t)>\liminf_{t\to +\infty}{f_2}'(t)=d_2,\]
we can write $^{\textsc{rl}}D^\alpha f\,(t)=({f_1}'(t)+f_3(t))+({f_2}'(t)-f_3(t))$, with $f_3\in\mathcal{C}_T(\mathbb{R}^+,\mathbb{R})$ and
\[d_1=\limsup_{t\to +\infty}{f_3}(t),\qquad d_2=\liminf_{t\to +\infty}{f_2}'(t).\]
This finishes the proof of this theorem.
\end{proof}
\end{theorem}

It seems fundamental to know when the fractional integral of order $\alpha\in(0,1)$ of a $T$-periodic function $f$ is a bounded function. Hereinafter we study this question.

\begin{theorem}\label{T:2.7}
Let $f\in\mathcal{C}_T(\mathbb{R},\mathbb{R})$ and $\alpha\in(0,1)$. Then, $I^\alpha f$ is a bounded function if, and only if,
\begin{equation}\label{eq:9}
\int_0^T f(t)\,dt=0.
\end{equation}
\begin{proof}
Let us suppose at the beginning that $I^\alpha f$ is a bounded function.

Since $f$ is a $T$-periodic function, using properties of Laplace transform, we obtain that
\[s\,\mathcal{L}\left[I^\alpha f\right](s)=\dfrac{s}{s^\alpha}\mathcal{L}[f](s),\]
where $\mathcal{L}\left[I^\alpha f\right]$ and $\mathcal{L}[f]$ are the Laplace transform of $I^\alpha f$ and $f$ respectively. Therefore, using Theorem \ref{T:1} and \cite[Theorem 1]{gluskin}, we have that
\begin{equation*}
\lim_{s\to 0^+}s\,\mathcal{L}\left(I^\alpha f\right)(s)=\dfrac{1}{T}\int_0^Tf(t)\,dt=\lim_{s\to 0^+}\dfrac{1}{s^\alpha}\dfrac{1}{T}\int_0^Tf(t)\,dt.
\end{equation*}
Hence, it is clear that if $I^\alpha f$ is bounded, then (\ref{eq:9}) holds.\medskip

Let us suppose now that (\ref{eq:9}) holds and we will prove that $I^\alpha f$ is, in such case, a bounded function. According to (\ref{eq:4}), we have that
\begin{equation}\label{eq:12}
I^\alpha f\,(t)={}^{\textsc{w}}I^\alpha f\,(t)-\lim_{\substack{n\to \infty\\n\in\mathbb{N}}}\int_{-nT}^0 (t-s)^{\alpha-1}f(s)\,ds.
\end{equation}
Since the Weyl integral of a $T$-periodic function which satisfies (\ref{eq:9}) is also a $T$-periodic function (see \cite[p. 348]{samko1993fractional}), we have that the first summand of the right hand side of (\ref{eq:12}) is bounded. Moreover, for the second summand, we have for each $n\in\mathbb{N}$, that
\begin{align}\label{eq:13}
\nonumber
\int_{-nT}^0(t&-s)^{\alpha-1}f(s)\,ds=\sum_{j=1}^n\int_{-jT}^{(1-j)T}(t-s)^{\alpha-1}f(s)\,ds\\
\nonumber
&=\sum_{j=1}^n\int_{0}^{T}(t-r+jT)^{\alpha-1}f(r-jT)\,dr
=\sum_{j=1}^n\int_{0}^{T}(t-r+jT)^{\alpha-1}f(r)\,dr\\
\nonumber
&=\sum_{j=1}^n\int_{0}^{T}(t-r+jT)^{\alpha-1}f^+(r)\,dr-\sum_{j=1}^n\int_{0}^{T}(t-r+jT)^{\alpha-1}f^-(r)\,dr\\
&\le c\left(\sum_{j=1}^n (t+(j-1)T)^{\alpha-1}-\sum_{j=1}^n(t-jT)^{\alpha-1}\right)=c\,t^{\alpha-1},
\end{align}
where
\begin{equation}\label{eq:c}
c=\int_0^Tf^+(t)\,dt=\int_0^Tf^-(t)\,dt>0.
\end{equation}

Analogously, for each $n\in\mathbb{N}$, we can obtain
\begin{equation}\label{eq:14}
\int_{-nT}^0(t-s)^{\alpha-1}f(s)\,ds\ge -c\,t^{\alpha-1}.
\end{equation}
So, from (\ref{eq:13}) and (\ref{eq:14}), we deduce that
\[-c\,t^{\alpha-1}\le \int_{-nT}^0(t-s)^{\alpha-1}f(s)\,ds\le c\,t^{\alpha-1}\]
and then, by Squezze Theorem \cite[Theorem 3.2.7]{bartle},
\[\lim_{t\to\infty} \lim_{\substack{n\to\infty\\ n\in\mathbb{N} }}\int_{-nT}^0(t-s)^{\alpha-1}f(s)\,ds=0.\]

Hence, we have proved that $I^\alpha f$ is a bounded function. This finishes the proof.
\end{proof}
\end{theorem}

\begin{remark}
Observe that in the second part of the proof of Theorem \ref{T:2.7}, we have proved that if $f\in\mathcal{C}_T(\mathbb{R},\mathbb{R})$, $\alpha\in(0,1)$ and (\ref{eq:9}) holds, then $I^\alpha f$ is an asymptotically $T$-periodic function. But note also that in the first part of the proof of Theorem \ref{T:2.7}, we have used the thesis of Theorem \ref{T:1}.
\end{remark}

\begin{figure}[h!]
\centering
\includegraphics[width=0.35\textwidth]{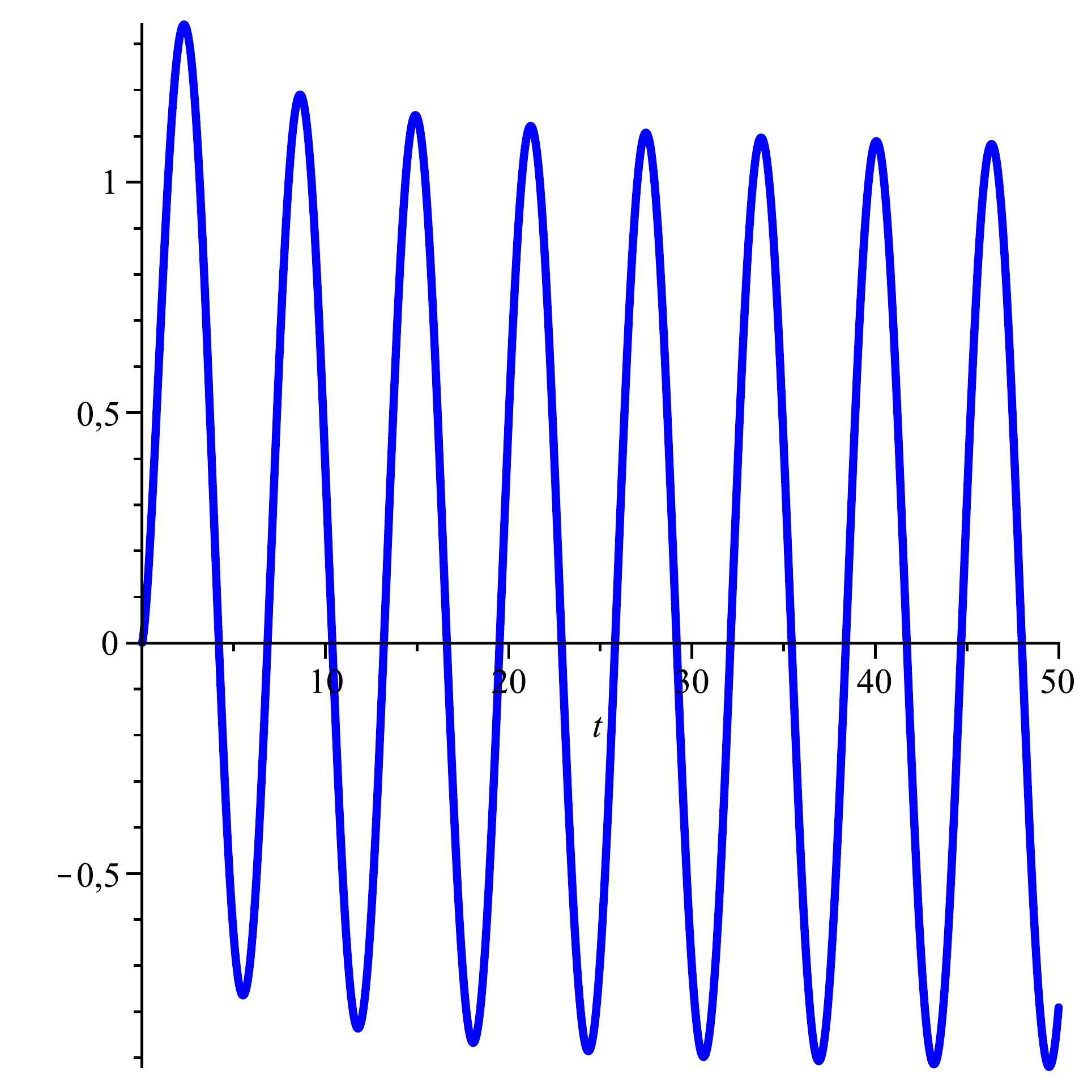}\hspace*{0.2cm}
\includegraphics[width=0.35\textwidth]{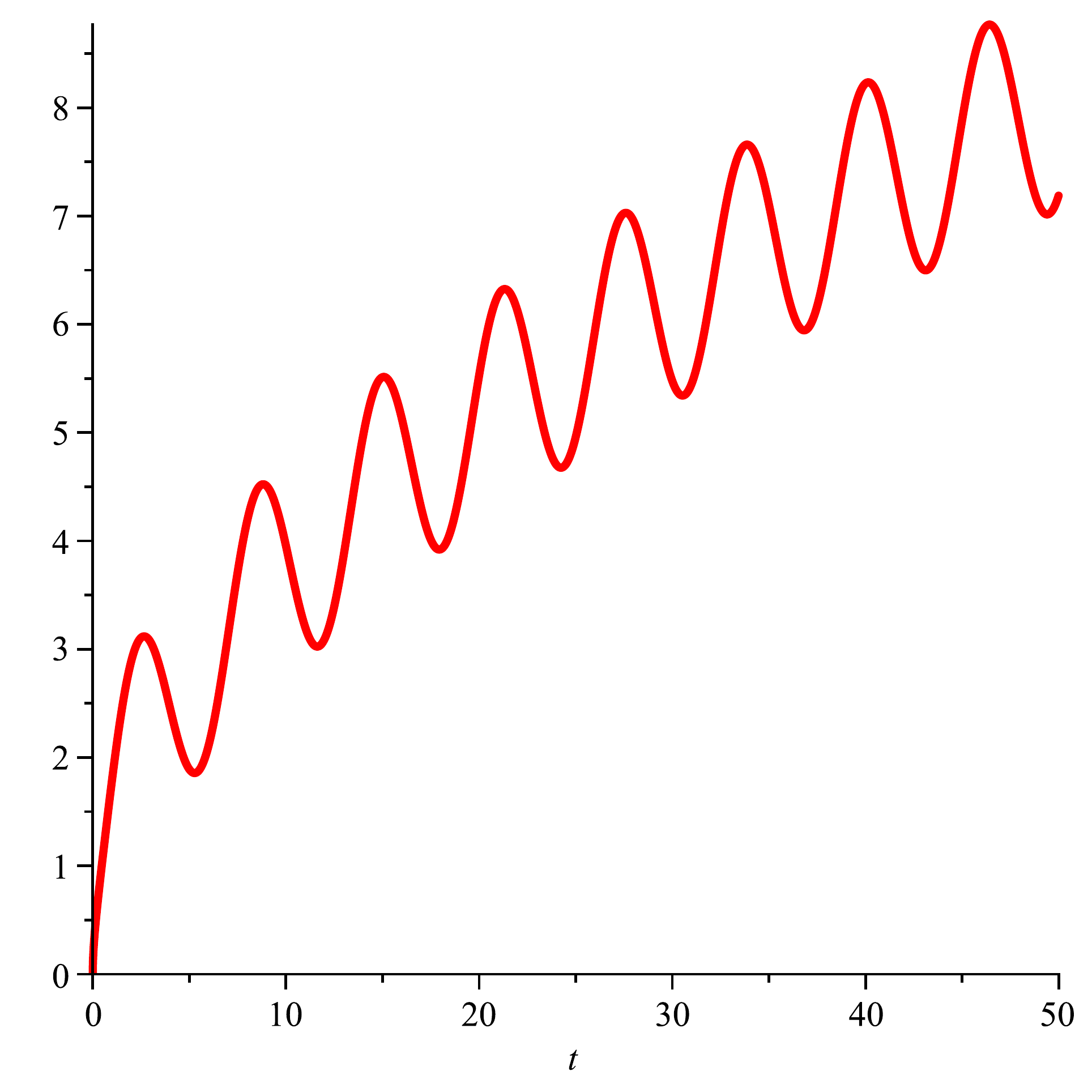}
\caption{Graph of $I^{1/2} \sin\,(t)$ (blue) and $I^{1/2}\left(\sin (r)+1\right)\,(t)$ (red) for $t\in[0,50]$.}\label{F:1}
\end{figure}

\begin{remark}
In \cite{MR2971825}, authors try to prove the thesis of Theorem \ref{T:2.7}. They use Terminal-Value Theorem of Laplace transform, but the proof is not correct because, using their notation, they suppose that exists $\lim_{t\to \infty}u(t)$ and this is not, in general, true. For example, for $\alpha\in(0,1)$, $\lim_{t\to\infty}I^\alpha\sin\,(t)$ does not exist (see Figure \ref{F:1}).
\end{remark}

\begin{remark}
For $\alpha>1$ we cannot say anything. For example, for $\alpha\in(1,2)$ we have,
\[I^\alpha \sin\,(t)=I^{\alpha-1}I^1\sin(t)=I^{\alpha-1}\left(-\cos s+1\right)\,(t),\]
which is an unbounded function. However,
\[I^\alpha \cos\,(t)=I^{\alpha-1}I^1\cos(t)=I^{\alpha-1}\sin\,(t).\]
So, $I^\alpha \cos$ is a bounded function.
\end{remark}

\begin{corollary}
Let $\alpha\in (0,1)$. If $f\in\mathcal{C}_T(\mathbb{R},\mathbb{R})$ is such that (\ref{eq:9}) does not hold, then
\[\lim_{t\to\infty}I^\alpha f\,(t)={\rm sgn}\left(\int_0^T f(t)\,dt\right)\infty.\]
\begin{proof}
Consider $\tilde{c}=c/T$, where $c$ is given by (\ref{eq:c}), and define the function $\tilde{f}$ as $\tilde{f}=f-\tilde{c}$. Now $\tilde{f}\in\mathcal{C}_T(\mathbb{R},\mathbb{R})$ and (\ref{eq:9}) holds for $\tilde{f}$. Since,
\[I^\alpha f\,(t)= I^\alpha\tilde{f}\,(t)+I^\alpha\tilde{c}\,(t)=I^\alpha\tilde{f}\,(t)+\dfrac{\tilde{c}}{\Gamma(1+\alpha)}\,t^\alpha,\]
the thesis of the corollary is now clear.
\end{proof}
\end{corollary}

We consider now the question of almost periodic functions.\medskip

By Proposition \ref{L:qp2} we have that, $AAP(\mathbb{R})\cap SAP_T(\mathbb{R})\subset AP_T(\mathbb{R})$. In addition, it is easy check that $AP_T(\mathbb{R})\subset SAP_T(\mathbb{R})$ holds (but $AP_T(\mathbb{R})\neq SAP_T(\mathbb{R})$, see \cite{henriquez1, nicola-pierri}) and it is trivial that $AP_T(\mathbb{R})\subset AAP(\mathbb{R})$. Hence, we have that
\[AAP(\mathbb{R})\cap SAP_T(\mathbb{R})= AP_T(\mathbb{R}).\]
But, by Proposition \ref{L:qp3}, we can deduce that
\begin{equation}\label{eq:15}
AP(\mathbb{R})\cap SAP_T(\mathbb{R})= P_T(\mathbb{R}),
\end{equation}
where $P_T(\mathbb{R})$ denotes the set of $T$-periodic functions.

Using these facts we can deduce the following results.
\begin{theorem}
If $f\in\mathcal{C}_T(\mathbb{R},\mathbb{R})$ and $I^\alpha f$, with $\alpha\in (0,1)$, is a bounded function then $I^\alpha f$ is not an almost periodic function.
\begin{proof}
Let us suppose that $I^\alpha f$ is an almost periodic function. In addition, by Theorem \ref{T:4.2}, $I^\alpha f$ is an S-asymptotically $T$-periodic function. Hence, by (\ref{eq:15}), $I^\alpha f$ is a $T$-periodic function, which contradicts Theorem \ref{T:noperiodic}.
\end{proof}
\end{theorem}

\begin{remark}
Reasoning in a similar way, we can obtain similar results for Riemann-Liouville and Caputo fractional derivatives of order $\alpha\in (0,1)$.
\end{remark}

\section{Example}\label{S:5}

Let, for each $t\in\mathbb{R}$, $f(t)=\sin(t)$ and let us consider the fractional integral of $f$ of order $\alpha>0$. Taking into account that
\[\sin(t)=\sum_{n=0}^{\infty} \frac{(-1)^n t^{2 n+1}}{(2 n+1)!}\]
and that for $\eta(t)=(-1)^{n} t^{2n+1}/(2n+1)!$, we have
\begin{equation*}
I^{\alpha}\eta\,(t)=\frac{1}{\Gamma(\alpha)}\int_{0}^{t}(t-s)^{\alpha-1} \frac{(-1)^n s^{2 n+1}}{(2 n+1)!} ds=\frac{(-1)^n t^{\alpha +2 n+1}}{\Gamma (\alpha +2 n+2)};
\end{equation*}
we finally deduce,
\begin{align*}
I^{\alpha}f\,(t)&=\frac{1}{\Gamma(\alpha)}\int_{0}^{t}(t-s)^{\alpha-1} f(s)ds=\frac{t^{\alpha +1} }{\Gamma (\alpha +2)}
\sum_{j=0}^{\infty} \frac{4^{-j} \left(-t^2\right)^j}{\left(\frac{\alpha }{2}+1\right)_j \left(\frac{\alpha+3 }{2}\right)_j}\\
&=\frac{t^{\alpha +1}}{\Gamma (\alpha +2)}\, _1F_2\left(1;\frac{\alpha }{2}+1,\frac{\alpha }{2}+\frac{3}{2};-\frac{t^2}{4}\right),
\end{align*}
assuming $t>0$, where $(A)_{j}=A\,(A+1)\cdots(A+j-1)$, $(A)_{0}:=1$ denotes the Pochhammer symbol.

For $\alpha=1/2$ the series expansion at infinity gives
\[I^{1/2}f\,(t) \sim \frac{\sqrt{\frac{1}{t}}}{\sqrt{\pi }}+\frac{\sin (t)-\cos (t)}{\sqrt{2}}\]
and therefore, we have that $I^{1/2}f$ is bounded (it has no poles).

Moreover, for $\alpha=3/2$ we obtain
\[I^{3/2}f\,(t) \sim \frac{2 \sqrt{t}}{\sqrt{\pi }}-\frac{\sin (t)+\cos (t)}{\sqrt{2}},\]
so $I^{3/2} f$ is unbounded.

Since for general $\alpha>0$ the first term at infinity is
\[\frac{t^{\alpha -3} \left(-(\alpha -3) \alpha +t^2-2\right)}{\Gamma (\alpha )}\]
we have that,  $I^{\alpha}f$ is bounded if and only if $\alpha \in (0,1]$.

\medskip

We can summarize that, for $\alpha \in (0,1)$ the function $I^{\alpha} \sin\,(t)$ is a bounded and S-asymp\-to\-ti\-cally $2\pi$-periodic function. Moreover, in these conditions, $I^\alpha \sin \,(t)$ is an asymptotically $2\pi$-periodic function, but it is not an almost periodic function.

\section{Concluding remarks}\label{S:7}

In this paper we have studied quasi-periodic properties of fractional derivatives and fractional integrals of periodic functions. Considering Riemann-Li\-ou\-vi\-lle and Caputo definitions, we have discussed conditions such as S-a\-symp\-to\-tic periodicity, asymptotic almost periodicity or asymptotic periodicity. We have also studied when the fractional integral and the fractional derivative of a $T$-periodic function is a bounded function.

Summarizing, we have obtained that, for $\alpha\in (0,1)$, the fractional integral $I^\alpha f$ of a $T$-periodic function is a bounded function if, and only if, the mean of $f$ is zero.
We have also provide conditions over $f$ which implie that $I^\alpha f$ is an S-asymptotically periodic or an asymptotically periodic function. These results concerning fractional integrals are useful to obtain similar results about fractional derivatives. We also like to note that these results could be helpful to a better understanding of some phenomena, such as for example those treated in \cite{brusselator,MR2863974}.

\section*{Funding}

The work of I. Area has been partially supported by the Ministerio de Econom\'{\i}a y Competi\-tividad of Spain under grant MTM2012--38794--C02--01, co-financed by the European Community fund FEDER. J. Losada and J.J. Nieto acknowledge partial financial support by the Ministerio de Econom\'{\i}a y Compe\-ti\-tividad of Spain under grants MTM2010--15314 and MTM2013--43014--P, Xunta de Galicia under grants R2014/002 and  Plan I2C ED481A-2015/272 co-financed by the European Community fund FEDER

\end{document}